\documentclass[11pt]{article}
\usepackage{graphicx}
\usepackage{amssymb}
\usepackage{color}
\usepackage[titletoc,title]{appendix}
\usepackage{amsmath}
\input{amssym.def}
\input{amssym.tex}

\newcommand{\I}{\mathop{\rm Im}\nolimits}

\begin{document}

\title{Vortex generated fluid flows in multiply connected domains}
\author{%
Anna Zemlyanova, Ian Manly, and Demond Handley\\
Department of Mathematics, Kansas State University, \\
138 Cardwell Hall, Manhattan KS 66506\\
Tel.: +1-785-532-6750, Fax: +1-785-532-0546}

\maketitle

\noindent

\sloppy

\begin{abstract}
A fluid flow in a multiply connected domain generated by an arbitrary number of point vortices is considered. A stream function for this flow is constructed as a limit of a certain functional sequence using the method of images. The convergence of this sequence is discussed, and the speed of convergence is determined explicitly. The presented formulas allow for the easy computation of the values of the stream function with arbitrary precision in the case of well-separated cylinders. The considered problem is important for applications such as eddy flows in the oceans. Moreover, since finding the stream function of the flow is essentially identical to finding the modified Green's function for Laplace's equation, the presented method can be applied to a more general class of applied problems which involve solving the Dirichlet problem for Laplace's equation.



\end{abstract}

\setcounter{equation}{0}

\section{Introduction}
A problem of fluid motion in the presence of vortices has important applications in geophysics, namely in the study of eddy flows in oceans. Ocean vortices may propagate large distances and are likely to encounter geographic obstacles such as islands, ocean ridges, and coastal lines. Vortex flows can be important vehicles for mass, momentum, heat, and salinity transfer in the oceans. Thus, the study of the vortex flows in multiply connected domains is important for accurate modeling and prediction of ocean flows. 

Motion of vortices in simply connected flow domains is relatively well-studied. The stream functions for these flows can be obtained by using the celebrated method of images \cite{MilneThomson1968} in combination with appropriate conformal mapping. The simplest example of application of the method of images is the study of a single vortex flow around one cylinder or an infinite straight wall. The resulting flow can be obtained by placing an image vortex with an opposite circulation at the symmetric point with respect to the cylinder or the straight wall. Reviews of the recent results on the vortex flows in simply connected domains are available in \cite{Arefetal2002, Newton2002}.

Scientific literature on vortex flows in multiply connected domains is considerably more limited. It is necessary to note the work by Johnson and McDonald \cite{JohnsonMcDonald2004} which is dedicated to the vortex flows in doubly connected domains. The solution is obtained by first conformally mapping the flow domain onto an annulus, then onto a periodically repeated rectangle in the complex plane, and exploiting the properties of elliptic theta functions. The vortex motion near walls with gaps is considered in \cite{JohnsonMcDonald2004b, JohnsonMcDonald2005}. Again, only simply and doubly connected domains are considered.

Vortex flows in the domains of arbitrary connectivity have been studied by Crowdy and Marshall in their multiple works \cite{CrowdyMarshall2005, CrowdyMarshall2005b, CrowdyMarshall2006}. The solutions in these papers have been obtained for multiply connected circular and slit domains in terms of the transcendental Schottky-Klein prime function \cite{Baker1995}. The numerical computation of the Schottky-Klein prime function is based upon computing an infinite product which does not converge in all cases. The convergence and the speed of convergence depends on the well-separatedness of the cylinders. Alternatively, the Schottky-Klein prime function can be computed by using power series approximations centered at the centers of the cylinders \cite{CrowdyMarshall2007b} in a similar way to the computation of the first-type Green's function for Laplace's equation in a circular domains \cite{Trefethen2005}.

In the present paper, a fluid flow generated by an arbitrary number of vortices around an arbitrary number of cylinders with specified circulation around each cylinder is studied. The stream function of the flow is obtained by the application of the method of images. According to the authors' knowledge the construction presented in this paper has not been attempted before. The main difficulty with applying the method of images to multiply connected flow domains is in the fact that the set of image vortices becomes infinite. This problem has been successfully overcome in the present paper. The stream function of the flow is obtained in terms of the limit of a certain functional sequence. The condition under which this sequence converges is investigated and depends on the mutual location and distance between the cylinders (so called well-separatedness of the cylinders). The speed of the convergence is investigated as well. In particular, it is established that the functional sequence converges with the speed of a geometric series. The presented solution is easy to implement numerically, and the results have been compared on the example of doubly connected domains to those obtained by using the method of elliptic functions in \cite{JohnsonMcDonald2004}.

Finally, it is necessary to note that finding the stream function for the vortex flow in question is essentially equivalent to finding the modified Green's function for a multiply-connected flow domain \cite{CrowdyMarshall2007a}. Hence, the presented technique can be applied to a much broader range of problems which can be reduced to solving the Dirichlet problem for Laplace's equation in multiply-connected domain. In particular, the applications of this method include such areas as electrostatics, potential theory, gravitation, numerical analysis, and approximation theory. Some of the alternative methods of construction of the Green's function using the theory of functional equations or the Schwarz-Christoffel mappings are presented in \cite{EmbreeTrefethen1999, MityushevRogosin2000}.

\setcounter{equation}{0}

\section{A vortex flow in a multiply connected domain}

\begin{figure}
\begin{center}
\includegraphics[width=0.5\textwidth]{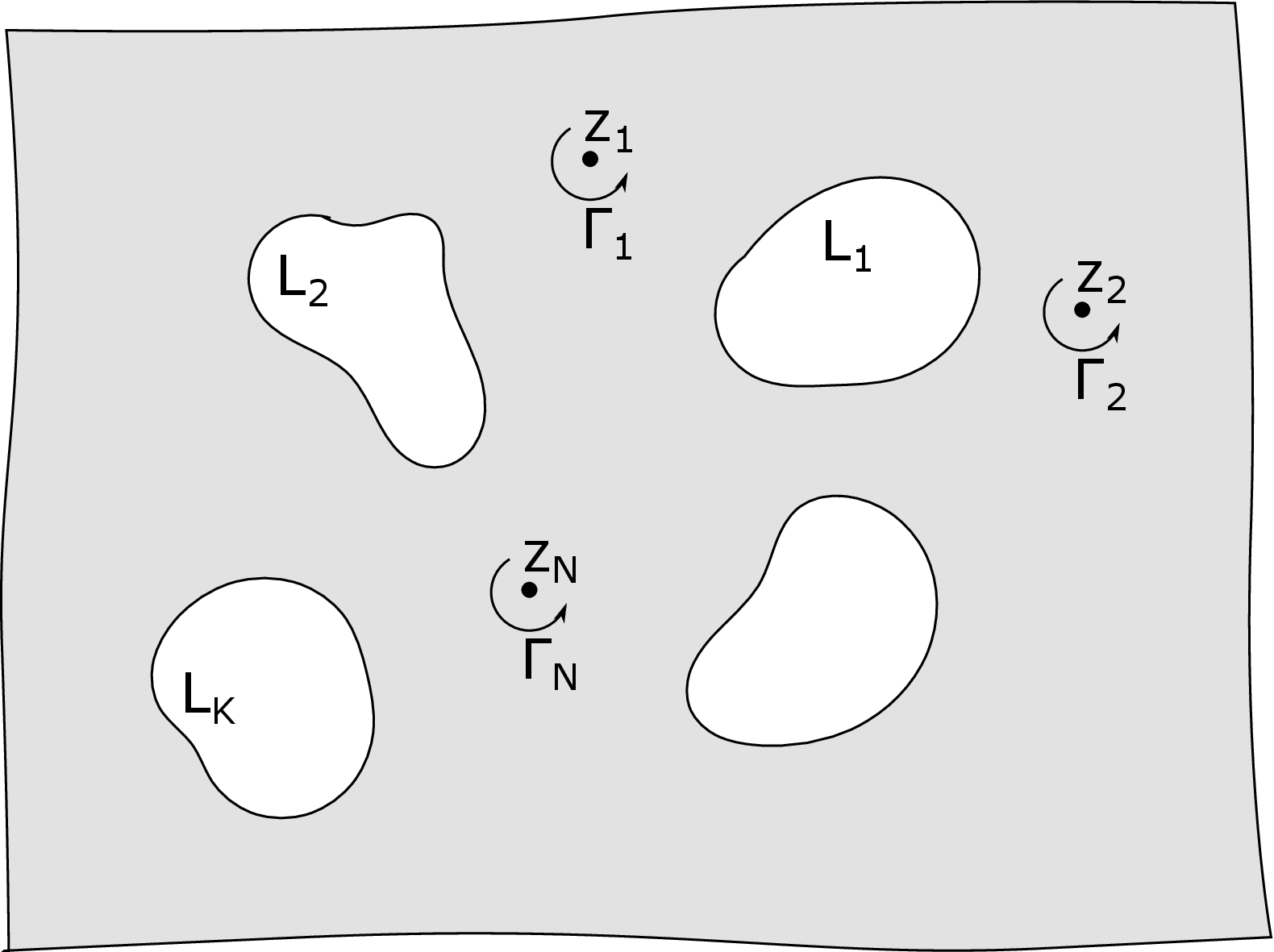} 
\end{center}
\caption{Vortex flow around $K$ islands in an unbounded domain.}
\label{fig1}
\end{figure}

Consider a flow of ideal fluid in an unbounded region $\tilde{D}$ exterior to $K$ islands $L_k$ of arbitrary smooth shape (fig. \ref{fig1}). The fluid flow in the domain $\tilde{D}$ is generated by $N$ point vortices located at the points $z_j$ with circulations $\Gamma_j$, $j=1,\ldots, N$. It is well known that a steady irrotational flow in the two-dimensional domain $\tilde{D}$ can be described by a complex potential $w(z)=\varphi(x,y)+i\psi(x,y)$ which is an analytic function in $\tilde{D}$ except for the points $z_j$ which are the singularities of the complex potential $w(z)$ of the logarithmic type:
\begin{equation}
w(z)=\frac{\Gamma_j}{2\pi i}\log (z-z_j)+\tilde{w}_j(z),
\label{2_1}
\end{equation}
where $\tilde{w}_j(z)$ is an analytic function in the neighborhood of the point $z_j$, $j=1,2,\ldots, N$.

The real $\varphi(x,y)$ and the imaginary $\psi(x,y)$ parts of the complex potential $w(z)$ are called correspondingly the velocity potential and the stream function of the flow. Both the velocity potential and the stream function are harmonic functions in $\tilde{D}\setminus \cup_{j=1}^N \{z_j\}$. Additionally, on the solid boundaries of the domain $\tilde{D}$ the stream function has to assume constant values:
\begin{equation}
\I w(z)=\mbox{Const},\,\,\, z\in L_j.
\label{2_2}
\end{equation}
The last condition from the physical viewpoint means that the solid boundaries are the streamlines of the flow.

Vortex trajectories in the region $\tilde{D}$ can be obtained by using the Kirchhoff-Routh function $H(x_1,y_1,\ldots,x_N,y_N)$. If $N$ vortices with circulations $\Gamma_j$, $j=1,\ldots,N$, are present in an incompressible fluid at the locations $(x_j(t),y_j(t))$ which depend on time $t$, then the trajectories of the vortices can be found from the following Hamiltonian equations \cite{Newton2002}:
$$
\Gamma_j\frac{dx_j}{dt}=\frac{\partial H}{\partial y_j},\,\,\,
\Gamma_j\frac{dy_j}{dt}=-\frac{\partial H}{\partial x_j}.
$$
The existence and the uniqueness of the Kirchhoff-Routh function $H(x_1,y_1,\ldots,x_N,y_N)$ has been established in \cite{Lin1941a}.
The relationship between the Kirchhoff-Routh function, the first-type Green's function and the complex potential of the flow has been described in detail in \cite{CrowdyMarshall2005, Lin1941a, Lin1941b, Newton2002}.

It should be noted that the complex potential and the stream function of the flow generated by several vortices in $\tilde{D}$ can be obtained by superposition of the flows generated by a single vortex in $\tilde{D}$. Thus, it is sufficient to consider the flow in the domain $\tilde{D}$ generated by only one vortex at the point $z_0$  with the unit circulation around this vortex, $\Gamma_0=1$. In this case the stream function $\psi(z)$ of the flow with zero circulations around each of the cylinders $L_j$ coincides with the modified Green's function \cite{CrowdyMarshall2005}.

Observe that the shape of the islands $L_k$ can be restricted to circular without loss of generality. By the generalization of the Riemann mapping theorem to multiply connected domains \cite{Nehari1952, Goluzin1969} there is a unique conformal mapping $f(z)$ of the $K$-connected domain $\tilde{D}$ onto some domain $D$ which is an exterior to $K$ circles in the extended complex plane $\overline{\mathbb{C}}$ (such a domain $D$ will be called from now on a circular domain) with the following expansion at infinity:
$$
f(z)=z+O\left(1/z\right).
$$
The circular domain $D$ is completely determined by the initial domain $\tilde{D}$ and the condition at infinity, and cannot be chosen arbitrarily. Since the conformal mapping $\omega=f(z)$ preserves the properties of the complex potential (\ref{2_1}), (\ref{2_2}), it can be assumed from now on that the flow of liquid is observed in the circular domain $D$. The corresponding flow in the original domain $\tilde{D}$ can be found then by taking the composition of the complex potential of the flow in $D$ with the conformal mapping $\omega=f(z)$ from  the original domain $\tilde{D}$ onto the circular domain $D$. While finding the exact analytic expression for the conformal mapping $\omega=f(z)$ is not feasible in most cases, very efficient numerical algorithms have been developed which allow to find this mapping approximately \cite{Delilloetal1999, Henrici1986}.

\begin{figure}
\begin{center}
\includegraphics[width=0.5\textwidth]{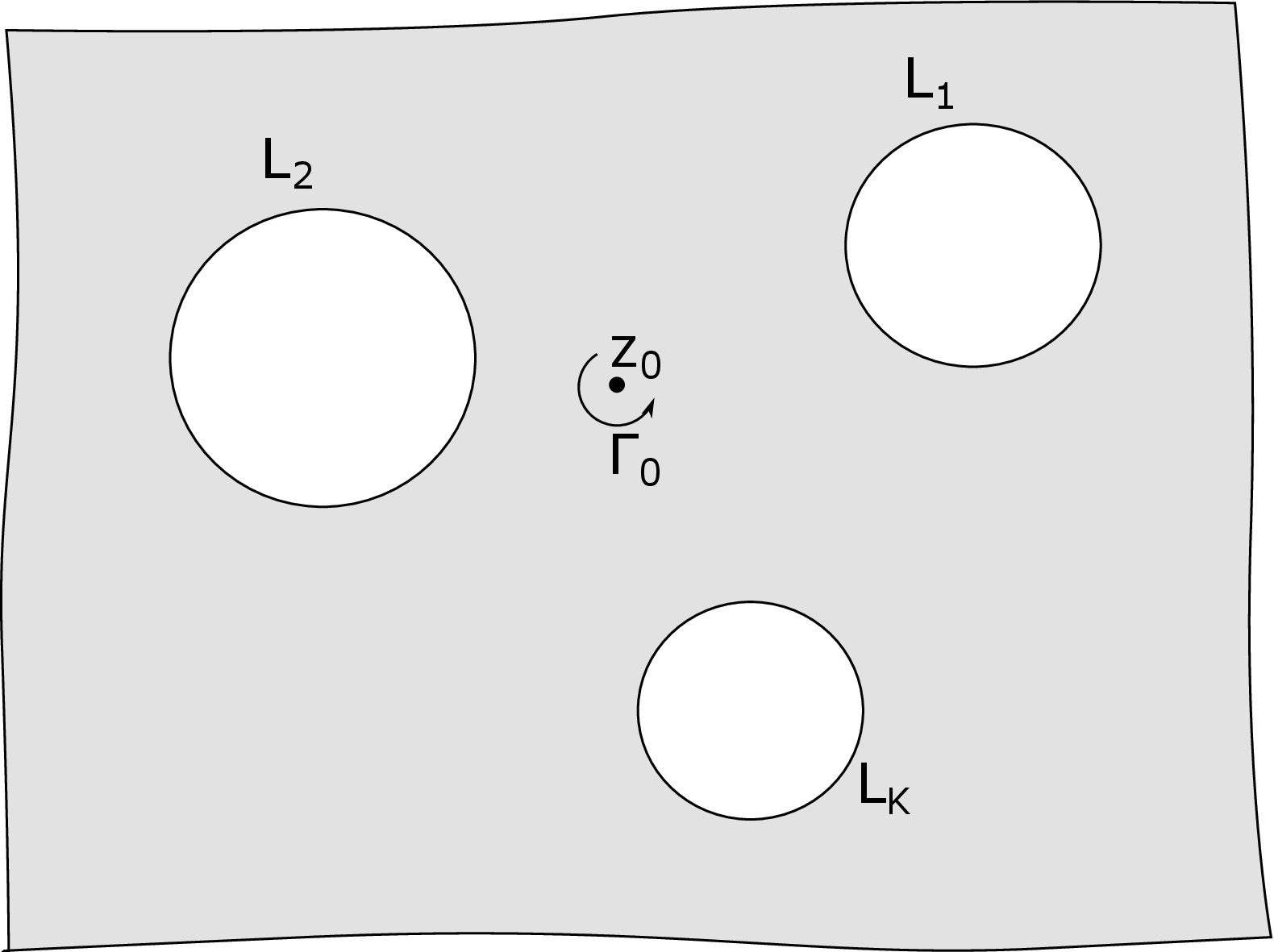} 
\end{center}
\caption{A single vortex in a $K$-connected circular domain.}
\label{fig2}
\end{figure}

From now on, consider the flow in the circular domain $D$ generated by a single vortex with a unit circulation $\Gamma_0=1$ located at the given point $z_0$ of the domain $D$ (fig. \ref{fig2}). Denote the stream function for this flow as $\psi^s(z,z_0)$. Then a flow in the circular domain $D$ generated by $N$ vortices located at the points $z_j$, $j=1,2,\ldots,N$, with circulations $\Gamma_j$, can be obtained by the superposition of the individual vortex flows for each of the points $z_j$:
\begin{equation}
\psi(z)=\sum_{j=1}^N \Gamma_j \psi^s(z,z_j).
\label{2_3}
\end{equation}

\setcounter{equation}{0}
\section{Method of images}

The stream function for the vortex flow shown on the fig. \ref{fig2} will be derived here by the application of the method of images. The main idea behind the method of images is to replace the original flow in the domain $D$ with impenetrable walls by a flow in the extended complex plane $\overline{\mathbb{C}}$ with additional ``image" vortices placed in the specially selected points of $\overline{\mathbb{C}}$ in such a way that the impenetrable walls of the original flow domain become the streamlines of the flow. 
\begin{figure}
\begin{center}
\includegraphics[width=0.3\textwidth]{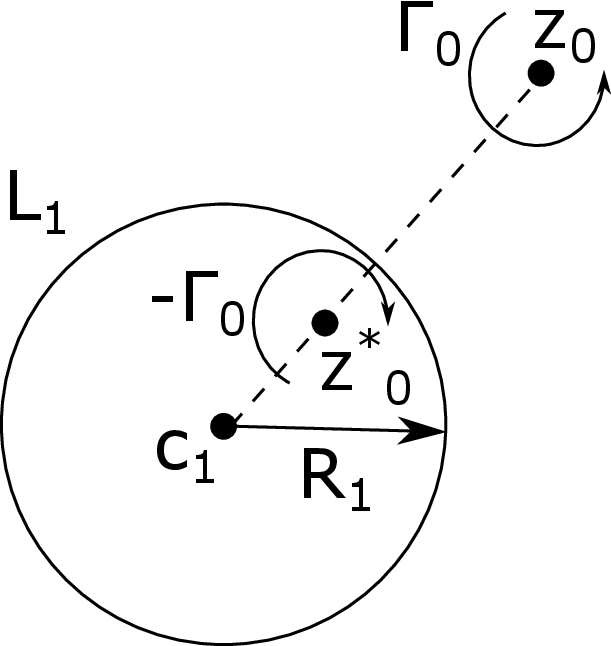} 
\end{center}
\caption{Method of images for one cylinder.}
\label{fig3}
\end{figure}

The method of images can be easily illustrated on the simple example of a vortex flow around one cylinder (fig. \ref{fig3}). Consider one vortex with a unit circulation $\Gamma_0=1$  which is located at the point $z_0$ of the complex plane, and one cylinder with impenetrable walls with a center $c_1$ and a radius $R_1$. To build a stream function for the flow around the cylinder, an image vortex needs to be placed at the point $z_0^*$ obtained by applying to $z_0$ the inversion map with respect to the circle $L_1:\, |z-c_1|=R_1$:
$$
T_1(z)=c_1+\frac{R_1^2}{\bar{z}-\bar{c_1}}
$$
Observe that  $z_0^*=T_1(z_0)$ and $T_1(z_0^*)=z_0$. The circulation at the image vortex is taken with the opposite sign $-\Gamma_0$ to that of the original vortex. Since the construction obtained in this way is symmetric with respect to the circle $L_1$, the circle $L_1$ becomes a streamline for the flow. This fact can be easily verified algebraically. The resulting stream function for the vortex flow around one cylinder has a form:
$$
\psi^s(z,z_0)=-\frac{1}{2\pi}\log |z-z_0|+\frac{1}{2\pi}\log|z-z_0^*|.
$$

The method of images in combination with a conformal mapping is easy to apply for the flows in simply-connected domains in $\overline{\mathbb{C}}$.
Using the method of images becomes more complicated in the case of multiple boundaries due to the fact that image vortices, in general, constitute an infinite set. Consider, in particular, the vortex flow depicted on the fig. \ref{fig2}. In this case there are $K$ inversion maps with respect to $K$ circles $L_j$:
\begin{equation}
T_j(z)=c_j+\frac{R_j^2}{\bar{z}-\bar{c}_j},\,\,\,j=1,\ldots,K.
\label{3_1}
\end{equation}
Again, observe that $T_j^2=I$, $j=1,2,\ldots,K$, where $I$ is an identity map.

\begin{figure}
\begin{center}
\includegraphics[width=0.5\textwidth]{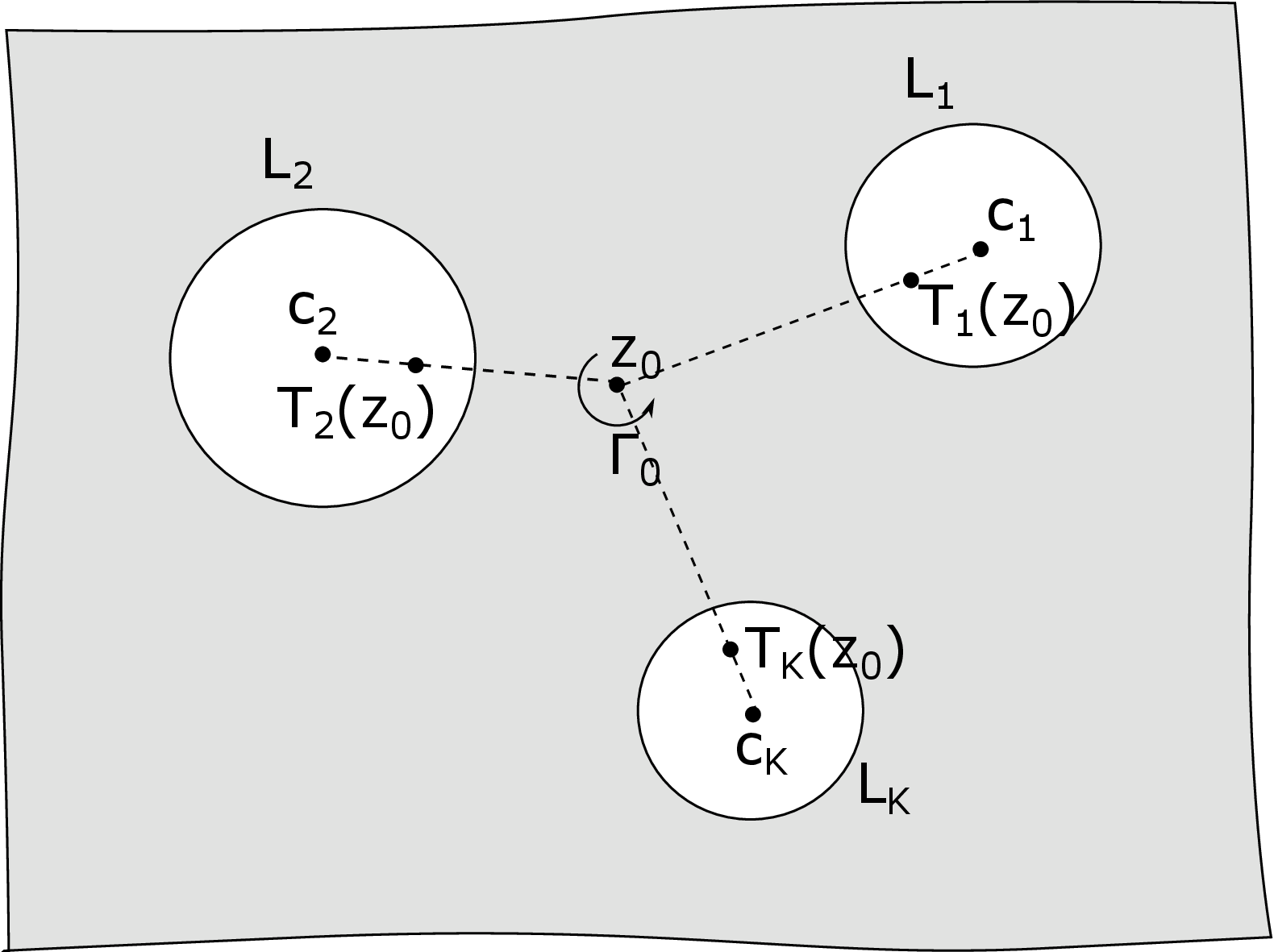} 
\end{center}
\caption{Level-1 symmetry points.}
\label{fig4}
\end{figure}

The goal of the method of images is to produce the set of ``image" vortices with respect to all of the rigid boundaries of the flow domain $D$. To obtain the image vortices, first, take the inversion maps of the point $z_0$ with respect to all $K$ circles. This produces level-1 symmetry points $T_j(z_0)$ shown on the fig. \ref{fig4}. 

However, unlike in the case of the single cylinder, this construction is not symmetric with respect to the circles $L_j$ since the points $T_j(z_0)$ do not have symmetric images inside the circles $L_k$, $k\neq j$. Thus, it is necessary to apply the inversion maps $T_j$ to the level-1 points which leads to the level-2 points $T_{i_1}T_{i_2}(z_0)$, $i_1\neq i_2$. Obviously, this process needs to be continued infinitely. Level-$N$ point can be written in the form $T_{i_1}T_{i_2}\ldots T_{i_N}(z_0)$ where $i_k\neq i_{k+1}$, $k=1,\ldots, N-1$. It is easy to count the symmetry points of each level obtained in this way: there are $K$ points of the level-1, $K(K-1)$ points of the level-2, $K(K-1)^{N-1}$ points of the level-$N$, and so on. Given that the original vortex has a unit circulation, the circulation at each of the symmetry points of the level-$N$ is equal to $(-1)^N$. Following the method of images, it is possible to formally write the stream function for the flow in the following way:
$$
\psi^s(z,z_0)=-\frac{1}{2\pi}\log|z-z_0|-\sum_{\zeta\in \mbox{level}\, 1}\frac{(-1)^1}{2\pi}\log|z-\zeta|-\ldots
$$
\begin{equation}
-\sum_{\zeta\in \mbox{level}\, N}\frac{(-1)^N}{2\pi}\log|z-\zeta|-\ldots
\label{3_2}
\end{equation}

Observe that without making additional assumptions about the summation, the formula (\ref{3_2}) does not make mathematical sense since it contains, in general, a divergent sum over the infinite set of the symmetry points $\zeta$. The sense in which the convergence of the formula (\ref{3_2}) is understood will be made precise soon.

Towards this purpose, consider the level-$N$ approximation to the sum (\ref{3_2}):
$$
\psi_N^s(z,z_0)=-\frac{1}{2\pi}\log|z-z_0|-\sum_{\zeta\in \mbox{level}\, 1}\frac{(-1)^1}{2\pi}\log|z-\zeta|-\ldots
$$
\begin{equation}
-\sum_{\zeta\in \mbox{level}\, N}\frac{(-1)^N}{2\pi}\log|z-\zeta|.
\label{3_3}
\end{equation}

Observe that even though the functions (\ref{3_3}) are harmonic in $D\setminus \{z_0\}$, and satisfy the condition (\ref{2_1}) at the point $z=z_0$, the condition (\ref{2_2}) will not be satisfied in the limit $N\to\infty$. Instead of $\psi_N^s(z,z_0)$ consider next the following function:
\begin{equation}
\psi^{s*}_N(z,z_0)=\frac{(K-1)\psi^s_N(z,z_0)+\psi^s_{N+1}(z,z_0)}{K}.
\label{3_4}
\end{equation} 

Our goal is to prove that the functional sequence $\psi^{s*}_N(z,z_0)$ converges to a harmonic function $\psi^s(z,z_0)$ in the domain $D\setminus \{z_0\}$. The resulting function $\psi^s(z,z_0)$ has a singularity of the logarithmic type (\ref{2_1}) at the point $z_0$ and is constant on the circles $L_j$, $j=1,\ldots, N$, which constitute the boundary of the domain $D$. In that case, the function $\psi^s(z,z_0)$ is the sought after stream function for the vortex flow in the circular domain $D$.

\setcounter{equation}{0}
\section{Convergence of the functional sequence $\{\psi_N^{s*}(z,z_0) \}_{N=1}^{\infty}$ in the domain $D\setminus \{z_0\}$}

Observe that the function $\psi^{s*}_N(z,z_0)$ can be rewritten in the following form:
\begin{equation}
-2\pi\psi^{s*}_N(z,z_0)=\frac{1}{K}\log|z-z_0|+
\label{4_1}
\end{equation}
$$
\frac{1}{K}\left((K-1)\log|z-z_0|+\sum_{\zeta\in\mbox{level}\, 1}(-1)^1\log|z-\zeta| \right)+
$$
$$
\frac{1}{K}\left((K-1)\sum_{\zeta\in\mbox{level}\, 1}(-1)^1\log|z-\zeta|+\sum_{\zeta\in\mbox{level}\, 2}(-1)^2\log|z-\zeta| \right)+\ldots
$$
$$
\frac{1}{K}\left((K-1)\sum_{\zeta\in\mbox{level}\, N}(-1)^N\log|z-\zeta|+\sum_{\zeta\in\mbox{level}\, N+1}(-1)^{N+1}\log|z-\zeta| \right).
$$

The last formula allows to split the expression for the function $\psi^{s*}_N(z,z_0)$ into the following ``layers":
\begin{equation}
(K-1)\sum_{\zeta\in\mbox{level}\, M}(-1)^M\log|z-\zeta|+\sum_{\zeta\in\mbox{level}\, M+1}(-1)^{M+1}\log|z-\zeta|.
\label{4_2}
\end{equation}

Since the number of the points of the level-$(M+1)$ is $K-1$ times larger than the number of the points of the level-$M$, then each layer (\ref{4_2}) contains the same number of terms corresponding to the levels $M$ and $M+1$. Any point of the level-$M$ can be written in the form ${\cal L}_M(z_0)=T_{i_1}T_{i_2}\ldots T_{i_M}(z_0)$ for some indices $i_1,i_2,\ldots,i_M\in \{1,\ldots, K \}$, $i_k\neq i_{k+1}$, $k=1,\ldots, (M-1)$. To each point of the level-$M$ correspond $K-1$ points of the level-$(M+1)$ which can be written in the form ${\cal L}_MT_{i_{M+1}}(z_0)$, $i_M\neq i_{M+1}$, $i_{M+1}\in \{1,\ldots, K \}$. Thus, each of the layers (\ref{4_2}) can be further split into the terms of the following type:
\begin{equation}
\log|z-{\cal L}_M(z_0)|-\log|z-{\cal L}_MT_{i_{M+1}}(z_0)|.
\label{4_3}
\end{equation}
Let us estimate each of the terms (\ref{4_3}). By applying the inequality 
$$
\log|1+x|\leq x\,\,\mbox{ for }\,\, x>-1,
$$ 
to the formula (\ref{4_3}), obtain:
$$
\log\left|\frac{z-{\cal L}_M(z_0)}{z-{\cal L}_MT_{i_{M+1}}(z_0)} \right|\leq \left|\frac{{\cal L}_MT_{i_{M+1}}(z_0)-{\cal L}_M(z_0)}{z-{\cal L}_MT_{i_{M+1}}(z_0)} \right|.
$$

The map ${\cal L}_M$ is a composition of $M$ inversions $T_j$, $j=i_1,i_2,\ldots,i_M \in \{1,\ldots, K \}$. Observe that at each inversion $T_j$, a couple of the symmetry points $z_1$, $z_2$ is mapped from the exterior of the given circle $L_j$ into the interior of this circle. Using the expression (\ref{3_1}), it is possible to write:
\begin{equation}
|T_j(z_1)-T_j(z_2)|=\frac{R_j^2|z_1-z_2|}{|z_1-c_j||z_2-c_j|}.
\label{4_4}
\end{equation}
Due to the procedure using which the symmetry points $z_1$, $z_2$ were generated, each of these points either lies in the interior of some other circle $L_k$ or is a point $z_0$. Simple geometric considerations show that:
$$
R_j/|z_1-c_j|<P_j,\,\,\,R_j/|z_2-c_j|<P_j,
$$
where
\begin{equation}
P_j=\max\left\{\frac{R_j}{|z_0-c_j|}, \,\,\,\frac{R_j}{|c_j-c_l|-R_l},\,l=1,\ldots,K,\,l\neq j \right\}.
\label{4_5}
\end{equation}
Then it is possible to conclude:
\begin{equation}
|{\cal L}_MT_{i_{M+1}}(z_0)-{\cal L}_M(z_0)|<P_{i_1}^2P_{i_2}^2\ldots P_{i_M}^2\frac{D}{R(z)}<P^{2M}\frac{D}{R(z)},
\label{4_6}
\end{equation}
where
\begin{equation}
P=\max\{P_1,\,P_2,\ldots,\,P_K \},
\label{4_7}
\end{equation}
$$
D=\max_{j}|T_j(z_0)-z_0|,
$$
$$
R(z)=\min_j\mbox{dist}(z,L_j),
$$
and $\mbox{dist}(z,L_j)$ denotes the shortest distance from the point $z$ to the circle $L_j$. Observe that for any compact set $K_0$ lying completely in the interior of the domain $D$ it is possible to find a number $R_0>0$ such that $R(z)>R_0$ for all the points $z\in K_0$.  Then the inequality (\ref{4_6}) can be rewritten as:
\begin{equation}
|{\cal L}_MT_{i_{M+1}}(z_0)-{\cal L}_M(z_0)|<P^{2M}\frac{D}{R_0}\,\,\, \mbox{for}\,\,\, \forall z\in K_0.
\label{4_8}
\end{equation}

Substituting (\ref{4_6}) and (\ref{4_8}) into the formula (\ref{4_1}), obtain that for any two positive integers $N_1$, $N_2$, such that $N_1<N_2$:
\begin{equation}
2\pi|\psi_{N_1}^{s*}(z,z_0)-\psi_{N_2}^{s*}(z,z_0)|<\frac{((K-1)P^2)^{N_1+1}D}{R(z)(1-(K-1)P^2)}, \,\,\,z\in D,
\label{4_9}
\end{equation}
or
\begin{equation}
2\pi|\psi_{N_1}^{s*}(z,z_0)-\psi_{N_2}^{s*}(z,z_0)|<\frac{((K-1)P^2)^{N_1+1}D}{R_0(1-(K-1)P^2)} \,\,\,\mbox{for}\,\,\, \forall z\in K_0,
\label{4_10}
\end{equation}
where $K_0\subset D$ is a compact set.

It follows from the inequality (\ref{4_9}) that the functional sequence $\psi_N^{s*}(z,z_0)$ is a Cauchy sequence in $D$ pointwise if the condition $(K-1)P^2<1$ holds. Hence, $\psi_N^{s*}(z,z_0)$ converges in $D$ pointwise to some function which we denote as $\psi^s(z,z_0)$. It follows from the inequality (\ref{4_10}) that the convergence is uniform on any compact set $K_0$. Hence, since the functions $\psi_N^{s*}(z,z_0)$ are harmonic in variable $z$ in $D\setminus \{z_0\}$ for all $N$ by construction, it follows that the limit function $\psi^s(z,z_0)$ is also harmonic in $D\setminus \{z_0\}$. The limit function $\psi^s(z,z_0)$ has a logarithmic singularity of the type (\ref{2_1}) at the point $z_0$ because all the functions $\psi_N^{s*}(z,z_0)$ have a singularity of this type at the point $z_0$.

The last property which needs to be proved is that the function $\psi^s(z,z_0)$ is constant on the circles $L_j$. Assume that $z\in L_j$ for some $j=1,2,\ldots,K$. Then from the properties of the inversion map $T_j$, it follows that $T_j(z)=z$. It can be seen that the points of the level-$(M+1)$ located inside any circle $L_j$ are obtained by taking an inversion map $T_j$ of the points of the level-$M$ located outside this circle:
$$
\{\zeta\in \mbox{level}\, (M+1),\,\,\zeta\in \mbox{int}\, L_j \}=T_j\{ \zeta\in \mbox{level}\, M,\,\,\zeta\notin \mbox{int}\, L_j \}.
$$

A simple algebraic computation then shows:
$$
\log|z-\zeta|-\log|z-T_j(\zeta)|=\log\left|\frac{\zeta-c_j}{R_j}\right|, \forall z\in L_j,
$$
where the right-hand side is constant on $L_j$.

Then we can rewrite the function $\psi_N^{s*}(z,z_0)$ as:
$$
-2\pi\psi_N^{s*}(z,z_0)=\log\left|\frac{z_0-c_j}{R_j}\right|+\sum_{\substack{\zeta\in\,\mbox{level}\, j,\\ j=1,\ldots,N-1,\\ \zeta\notin \,\mbox{int}\,L_j}}(-1)^j\log\left|\frac{\zeta-c_j}{R_j}\right|+
$$
\begin{equation}
\frac{1}{K}\sum_{\substack{\zeta\in\,\mbox{level}\, N,\\ \zeta\notin \mbox{int}\,L_j}}(-1)^N\log\left|\frac{\zeta-c_j}{R_j}\right|+
\label{4_11}
\end{equation}
$$
\frac{K-1}{K}\sum_{\substack{\zeta\in\,\mbox{level}\, N,\\ \zeta\notin \mbox{int}\,L_j}}(-1)^N\log\left|z-\zeta\right|+\frac{1}{K}\sum_{\substack{\zeta\in\,\mbox{level}\, (N+1),\\ \zeta\notin \mbox{int}\,L_j}}(-1)^{N+1}\log\left|z-\zeta\right|.
$$
Observe that the first three terms of the formula (\ref{4_11}) are independent of the point $z\in L_j$, and the last two terms can be estimated similarly to (\ref{4_9}):
$$
\left|\frac{K-1}{K}\sum_{\substack{\zeta\in\,\mbox{level}\, N,\\ \zeta\notin \mbox{int}\,L_j}}(-1)^N\log\left|z-\zeta\right|+\frac{1}{K}\sum_{\substack{\zeta\in\,\mbox{level}\, (N+1),\\ \zeta\notin \mbox{int}\,L_j}}(-1)^{N+1}\log\left|z-\zeta\right|\right|<
$$
$$
\frac{K-1}{K}\frac{((K-1)P^{2})^ND}{R_j(z)},
$$
where 
$$
R_j(z)={\min}_{\substack{k=1,\ldots,K,\\k\neq j}}\mbox{dist}(z,L_k).
$$

Thus, under assumption $(K-1)P^2<1$, it follows that the values of the sequence $\psi_N^{s*}(z,z_0)$ converge to constants on $L_j$ for each $j=1,\ldots,K$. It is possible to show that these constants are finite. To do so, we again split the remaining terms in (\ref{4_11}) into the ``layers":
$$
\frac{K-1}{K}\sum_{\substack{\zeta\in\,\mbox{level}\, M,\\ \zeta\notin \mbox{int}\,L_j}}(-1)^j\log\left|\frac{\zeta-c_j}{R_j}\right|+\frac{1}{K}\sum_{\substack{\zeta\in\,\mbox{level}\, (M+1),\\ \zeta\notin \mbox{int}\,L_j}}(-1)^{N+1}\log\left|\frac{\zeta-c_j}{R_j}\right|,
$$
which can be further split into the individual terms and estimated:
$$
\log|{\cal L}_M(z_0)-c_j|-\log|{\cal L}_MT_{i_{M+1}}(z_0)-c_j|<\frac{P^{2M}D}{\min_{\substack{k=1,\ldots,N,\\k\neq j}}(|c_k-c_j|-R_k)}.
$$
From the last estimate it is possible to conclude that the values of $\psi^s(z,z_0)$ are finite on each of the circles $L_j$ if the condition $(K-1)P^2<1$ holds.

It follows then that the limit function $\psi^s(z,z_0)=\lim_{N\to\infty}\psi_N^{s*}(z,z_0)$ satisfies all of the conditions imposed on the stream function for the considered vortex flow in the $K$-connected circular domain $D$.

\setcounter{equation}{0}
\section{Circulations around cylinders and at infinity}
The stream function $\psi^s(z,z_0)$ corresponds to the flow in the domain $D$ with a single vortex of a unit circulation $\Gamma_0=1$ located at the point $z_0$. The flow in the domain $D$ with $N$ vortices located at the points $z_j$, $j=1,\ldots,N$, with circulations $\Gamma_j$, can be easily obtained from the stream function $\psi^s(z,z_0)$ by superposition of the stream functions for the individual vortices using the formula (\ref{2_3}).

Consider the circulations around the cylinders $L_j$ which are prescribed by the stream function $\psi^s(z,z_0)$. The circulation around a closed contour $C$ in a fluid domain can be computed by the following formula:
$$
\Gamma_C=\oint_C {\bf u} \cdot d{\bf s},
$$
where $\bf u$ is a velocity and $d{\bf s}$ is an element along the contour.

Using this formula it is possible to obtain that the circulation around any cylinder $L_j$ is equal to $-1/K$ for the functions $\psi_N^{s*}(z,z_0)$ for all $N$. Thus, in the limit $N\to\infty$, the circulation of the flow defined by the stream function $\psi^s(z,z_0)$ is also equal to $-1/K$ around any cylinder $L_j$ for all $j=1,\ldots,N$. 

A vortex at infinity point $z_0=\infty$ of the domain $D$ can be introduced by using a similar procedure as for the vortex at a finite point $z_0$.  In particular, the flow with a vortex at infinity with a given circulation $\Gamma_{\infty}$ can be generated by the formulas:
$$
\psi_{N,\infty}(z)=-\sum_{\zeta\in \mbox{level}\, 1}\frac{(-1)^1\Gamma_{\infty}}{2\pi}\log|z-\zeta|-\ldots
$$
\begin{equation}
-\sum_{\zeta\in \mbox{level}\, N}\frac{(-1)^N\Gamma_{\infty}}{2\pi}\log|z-\zeta|,
\label{5_1}
\end{equation}
\begin{equation}
\psi^*_{N,\infty}(z)=\frac{(K-1)\psi_{N,\infty}(z)+\psi_{N+1,\infty}(z)}{K},
\label{5_2}
\end{equation} 
and letting $N\to \infty$. These formulas are analogous to the formulas (\ref{3_3}), (\ref{3_4}) with the only exception that the first ``generating term" for the vortex at the point $z_0=\infty$ is omitted. Observe, that the level-1 points in this case are the centers $c_j$ of the circles $L_j$. Similarly to the case of a finite point $z_0$, placing a vortex at the point $z_0=\infty$ with the circulation $\Gamma_{\infty}$ induces the circulations equal to $-\Gamma_{\infty}/K$ around each of the cylinders $L_j$, $j=1,\ldots,K$.

Finally, for some practical applications, it is important to prescribe the circulations around the cylinders $L_j$, $j=1,\ldots,N$. This can be done by placing additional vortices with circulations $\Gamma_j^c$ at the centers $c_j$ of the cylinders. The formulas (\ref{3_3}), (\ref{3_4}) then become:
$$
\psi_{N,j}(z)=-\frac{\Gamma_j^c}{2\pi}\log|z-c_j|-\sum_{\substack{\zeta\in \mbox{level}\, 1,\\\zeta\neq\infty}}\frac{(-1)^1\Gamma_j^c}{2\pi}\log|z-\zeta|-\ldots
$$
\begin{equation}
-\sum_{\zeta\in \mbox{level}\, N}\frac{(-1)^N\Gamma^c_j}{2\pi}\log|z-\zeta|.
\label{5_3}
\end{equation}
\begin{equation}
\psi^*_{N,j}(z)=\frac{(K-1)\psi_{N,j}(z)+\psi_{N+1,j}(z)}{K}.
\label{5_4}
\end{equation}
In this case $\zeta=\infty$ is a level-1 point, and the corresponding term must be omitted in the formula (\ref{5_3}). The level-2 points contain all the centers $c_k$, $k\neq j$, which are the symmetry points of the infinity point with respect to the cylinders $L_k$, $k\neq j$. Again, placing the vortex with a circulation $\Gamma_j^c$ at the point $c_j$ induces a circulation $\Gamma^c_j(2K-1)/K$ around the cylinder $L_j$, additional circulations of $-\Gamma^c_j/K$ around all other cylinders $L_k$, $k\neq j$, and a circulation $-\Gamma^c_j$ at infinity.

Finally, combining the stream functions for individual vortices, vortex at infinity, and vortices at the centers of the cylinders $L_j$ one can obtain a flow in the domain with $N$ vortices with any prescribed circulations around each of the cylinders $L_j$, $j=1,\ldots,N$. If the desired circulations around each of the cylinders $L_j$ are equal to $\gamma_j$, this leads to the following system of linear algebraic equations with respect to the unknowns $\Gamma^c_j$:
\begin{equation}
\left[
\begin{array}{cccc}
\frac{2K-1}{K} & -\frac{1}{K} & \vdots & -\frac{1}{K}\\
-\frac{1}{K} & \frac{2K-1}{K} & \vdots & -\frac{1}{K}\\
\ldots & \ldots & \ddots & \ldots\\
-\frac{1}{K} & -\frac{1}{K} & \vdots & \frac{2K-1}{K}
\end{array}
\right]
\left[
\begin{array}{c}
\Gamma^c_1\\
\Gamma^c_2\\
\ldots\\
\Gamma^c_K
\end{array}
\right]=
\left[
\begin{array}{c}
\gamma_1+\frac{1}{K}\sum_{l=1}^N \Gamma_l+\frac{\Gamma_{\infty}}{K}\\
\gamma_2+\frac{1}{K}\sum_{l=1}^N \Gamma_l+\frac{\Gamma_{\infty}}{K}\\
\ldots\\
\gamma_K+\frac{1}{K}\sum_{l=1}^N \Gamma_l+\frac{\Gamma_{\infty}}{K}
\end{array}
\right]
\label{5_5}
\end{equation}

The system has a diagonally dominant matrix and, hence, is uniquelly solvable for any right-hand side. Finally, the circulation at infinity in this case will be equal to $-\sum_{l=1}^K\gamma_l-\sum_{l=1}^N\Gamma_l$.

\setcounter{equation}{0}
\section{Set of the symmetry points and its limit set}

Consider the set of all the symmetry points corresponding to the point $z=z_0$. These points can be described by the formula $T_{i_1}T_{i_2}\ldots T_{i_M}(z_0)$ for some $M\geq 1$ and for some set of indices $i_k\in \{1,\ldots, K \}$, $i_k\neq i_{k+1}$, $k=1,\ldots, (M-1)$. Let us study first the case of a doubly-connected domain with only two cylinders $L_1$, $L_2$ present. Observe that in this case all of the symmetry points can be described by four sequences:
$$
a_j=(T_1T_2)^{j-1}T_1(z_0),\,\,b_j=(T_1T_2)^j(z_0),
$$
$$
c_j=(T_2T_1)^{j-1}T_2(z_0),\,\,d_j=(T_2T_1)^j(z_0),\,\,\,j=1,2,\ldots.
$$

The points of the sequences $\{a_j\}_{j=1}^{\infty}$ and $\{b_j\}_{j=1}^{\infty}$ lie in the interior of the circle $L_1$, while the points of the sequences $\{c_j\}_{j=1}^{\infty}$ and $\{d_j\}_{j=1}^{\infty}$ lie in the interior of the circle $L_2$. It is easy to show that all four sequences are convergent, and the limit points are the fixed points of the mappings $T_1T_2(z)$ and $T_2T_1(z)$. In particular,
$$
a_j\to z^{\star}_1,\,\,b_j\to z^{\star}_1,\,\,c_j\to z^{\star}_2,\,\,d_j\to z^{\star}_2,\,\,\mbox{as}\,\, j\to\infty,
$$
where
$$
T_1T_2(z^{\star}_1)=z^{\star}_1,\,\,z^{\star}_1\in \,\mbox{int}\, L_1,\,\,\,T_2T_1(z^{\star}_2)=z^{\star}_2,\,\,z^{\star}_2\in\,\mbox{int}\, L_2.
$$

\begin{figure}
\begin{center}
\includegraphics[width=0.5\textwidth]{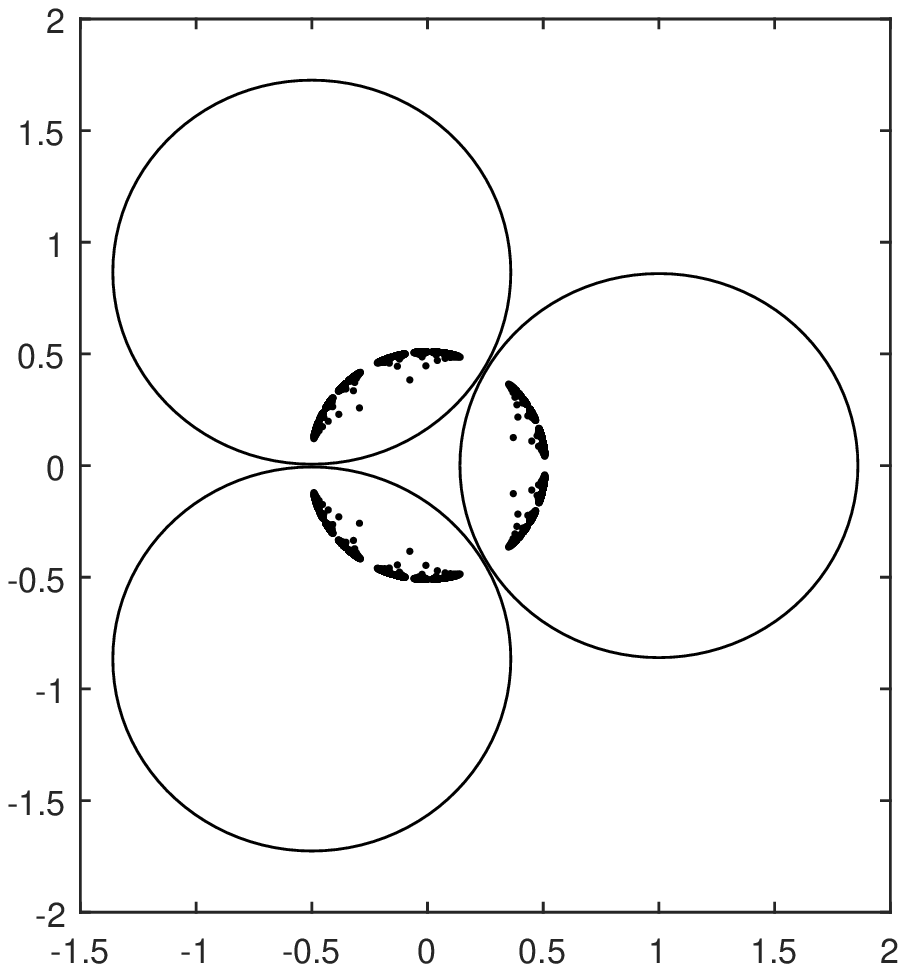} \,\,\,
\includegraphics[width=0.8\textwidth]{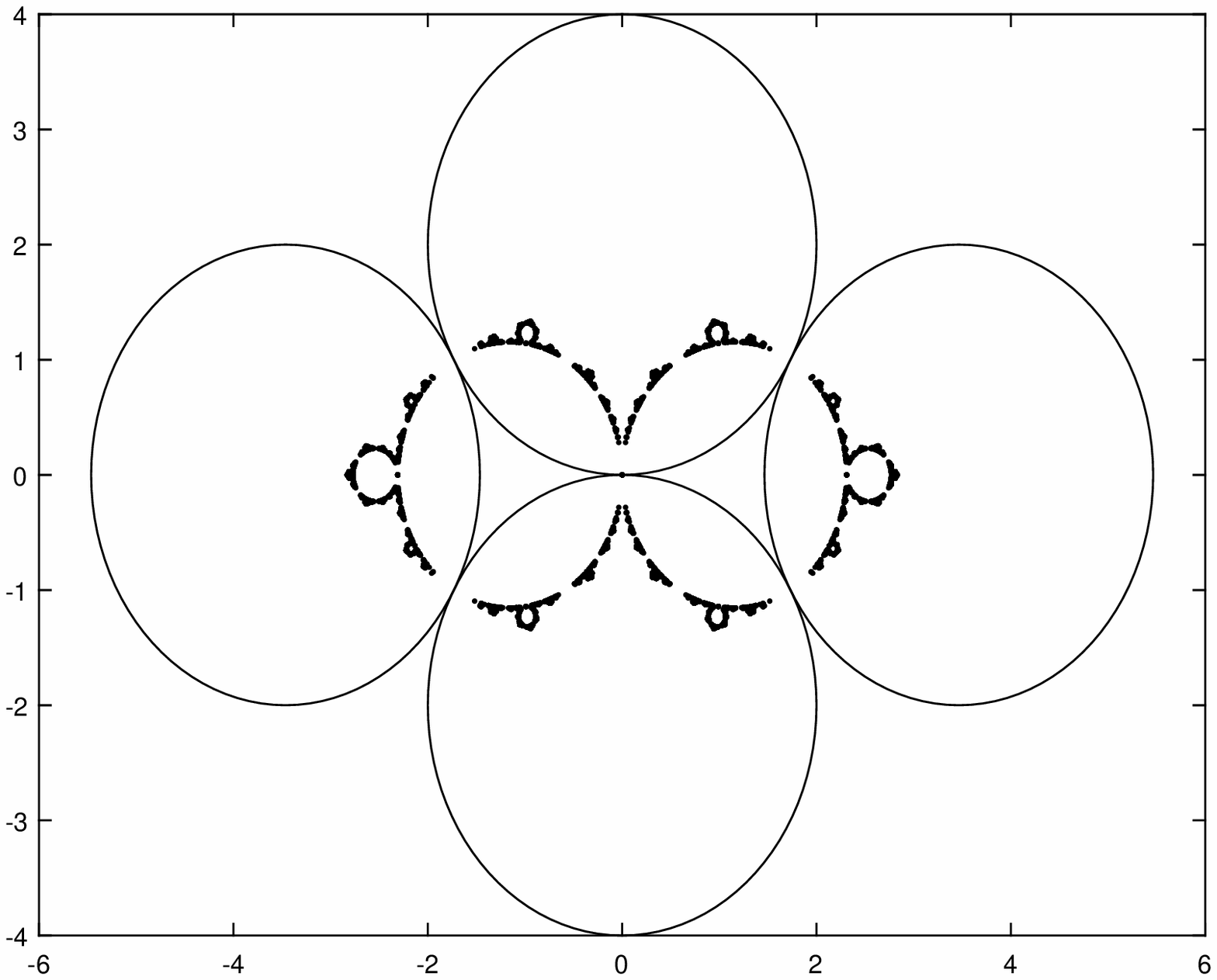} 
\end{center}
\caption{Sets of the symmetry points for touching circles.}
\label{fig5}
\end{figure}
The situation becomes more complicated for domains of connectivity higher than two. Observe that the set of symmetry points will necessary be self-similar. This follows from the fact that the level-$(N+1)$ points are obtained from the level-$N$ points by applying one of the symmetry maps $T_j$, $j=1,\ldots,K$. Observe also that in limiting cases, when the circles touch, the limiting  set of the symmetry points can become a circle or even a fractal (fig. \ref{fig5}). The fig. \ref{fig5} are plotted with symmetry points up to the level $10$. More information about the limiting sets of the symmetry maps and M\"obius maps can be found in \cite{Mumford2015}.

\setcounter{equation}{0}
\section{Numerical results}
\begin{figure}
\begin{center}
\includegraphics[width=0.45\textwidth]{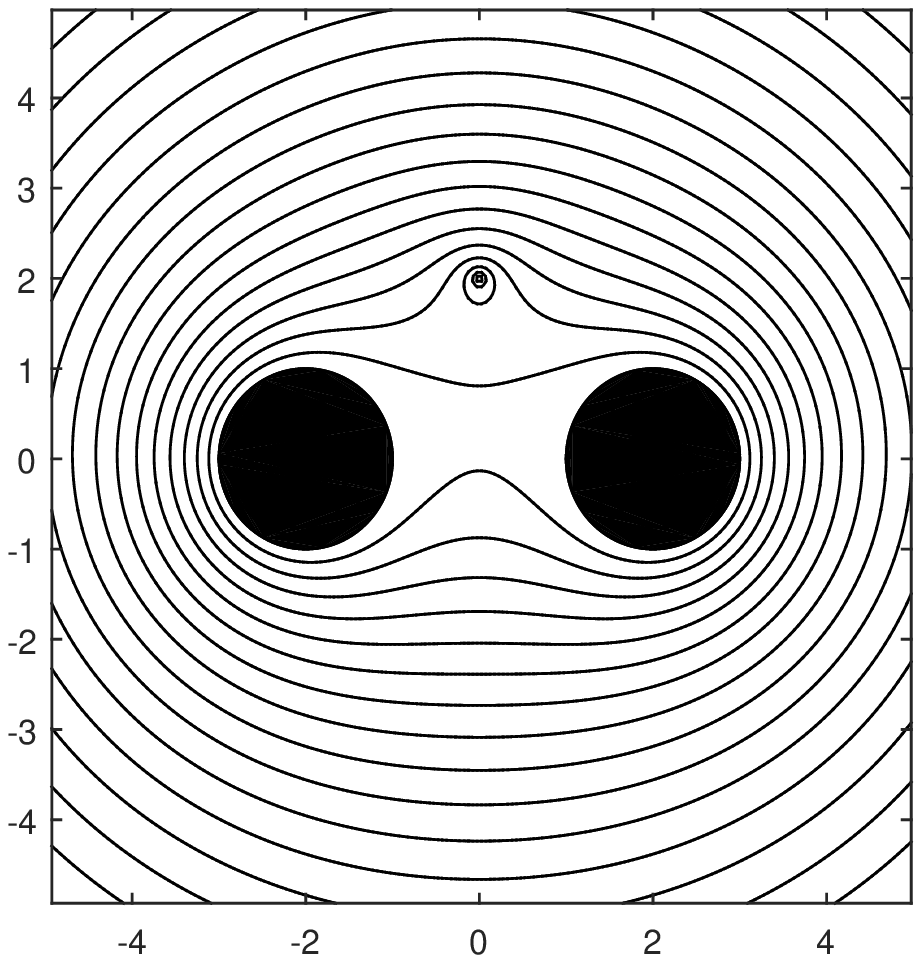}
\end{center}
\caption{Vortex flow around two cylinders with zero circulations on the boundaries of the cylinders.}
\label{fig6}
\end{figure}
Consider a vortex flow around two circular cylinders. The resulting flow domain in this case is doubly connected. Using a conformal mapping the exterior of two circles can be first mapped onto a concentric circular ring, and then onto a rectangle periodically repeated throughout the whole complex plane \cite{JohnsonMcDonald2004}. Hence, the solution can be furnished in terms of elliptic functions, namely, the elliptic theta function $\theta_1(\zeta,q)$. Numerical comparison of the results obtained by the methods used in \cite{JohnsonMcDonald2004} with the results of the current paper is given in the table \ref{tab1}. The computations in the table are made for the stream function of the vortex flow around two cylinders with the parameters $c_1=0$, $R_1=1$, $c_2=3$, $R_2=0.5$, $z_0=2i$, $\Gamma_0=1$ and zero circulations on both cylinders which is equivalent to placing a vortex with a circulation $\Gamma_{\infty}=-\Gamma_0$ at the infinity point of the plane.
\begin{center}
\begin{table}[!h]
\begin{tabular}{ |c|c|c| } 
\hline
 Point $z$ & Johnson and McDonald & Current paper \\ 
\hline
 $-3.5-3.5i$ & $-0.174608512540543$ & $-0.174608618004631$ \\ 
 $0.5-1.5i$ & $-0.047561219605849$ & $-0.047561611378318$ \\ 
 $2.5+3.5i$ & $-0.073398543207433$ & $-0.073398504917044$ \\
 $1.5+0.5i$ & $-0.020268684918721$ & $-0.020268607383453$ \\
 \hline
\end{tabular}
\label{tab1}
\caption{Comparison of the values of the stream function for a doubly connected domain.}
\end{table}
\end{center}
Observe that for the doubly connected flow domain the condition $(K-1)P^2<1$ is always satisfied and, hence, the method presented in this paper always converges irrespectively of the mutual location and the size of the cylinders and vortices.

\begin{figure}
\begin{center}
\includegraphics[width=0.45\textwidth]{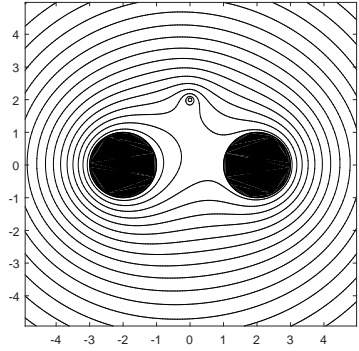}
\end{center}
\caption{Vortex flow around two cylinders with the circulation equal to $1/2$ on the left cylinder and the circulation equal to $-1/2$ on the right cylinder.}
\label{fig7}
\end{figure}
\begin{figure}
\begin{center}
\includegraphics[width=0.45\textwidth]{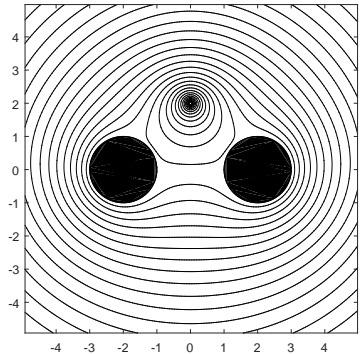}
\end{center}
\caption{Vortex flow around two cylinders with the circulations equal to $1/2$ on the boundaries of the cylinders.}
\label{fig8}
\end{figure}
The graphs of the instantaneous streamlines for the vortex flow around two cylinders are plotted on the figs. \ref{fig6}, \ref{fig7}, \ref{fig8}. The single vortex is located at the point $z_0=2i$ and has the circulation $\Gamma_0=1$, the cylinders have the centers $c_1=-2$, $c_2=2$, and the radii $R_1=R_2=1$. The circulations around both cylinders on the fig. \ref{fig6} are equal to zero, on the fig. \ref{fig7} the circulation around the left cylinder is equal to $1/2$ and around the right cylinder to $-1/2$, and on the fig. \ref{fig8} the circulations around both cylinders are equal to $1/2$.
\begin{figure}
\begin{center}
\includegraphics[width=0.45\textwidth]{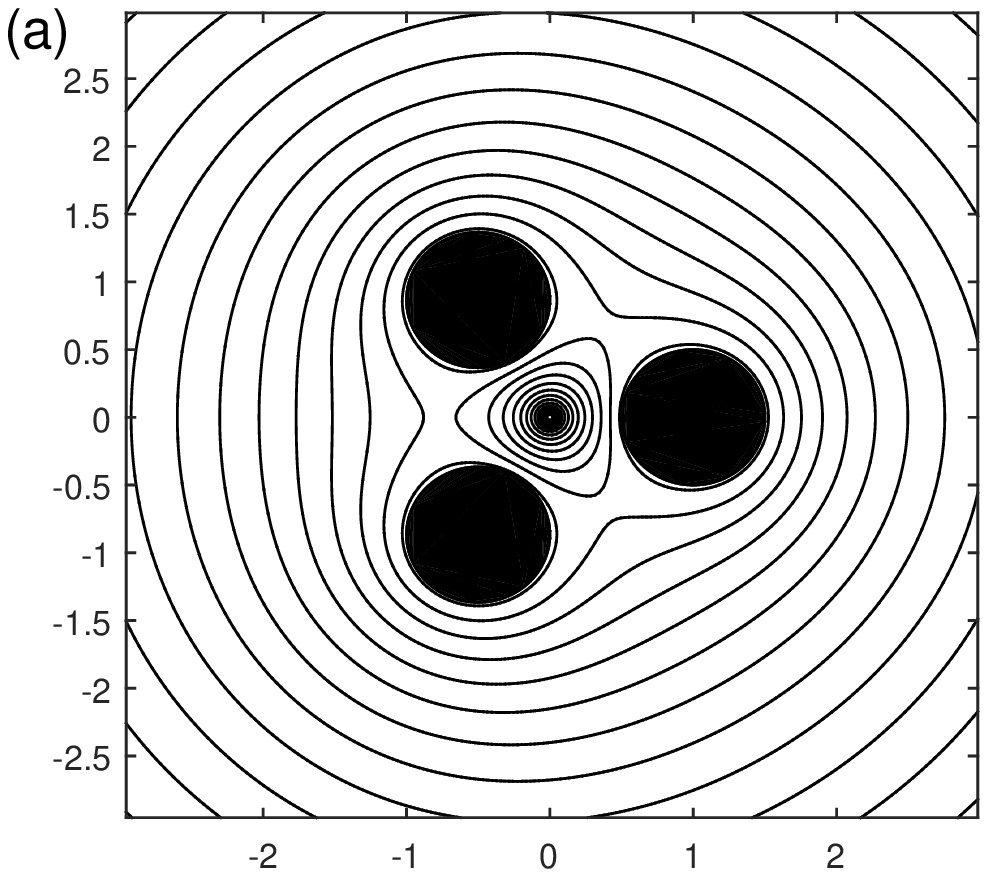}\,\,
\includegraphics[width=0.45\textwidth]{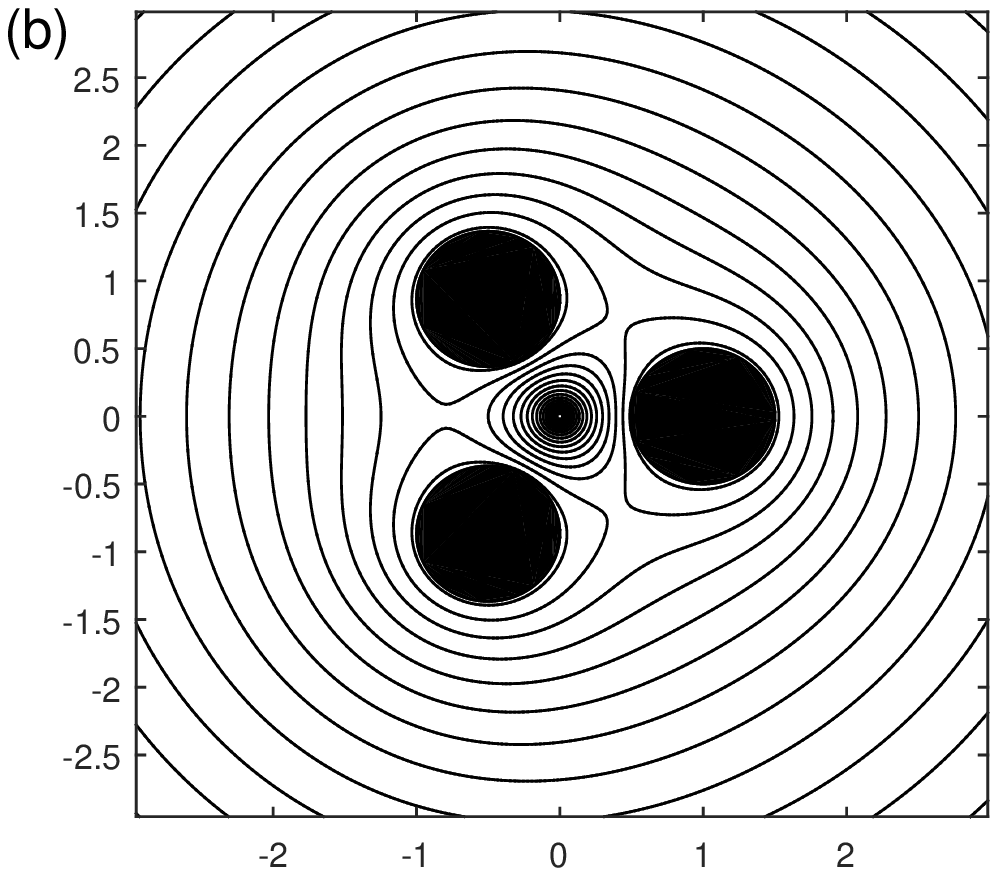}
\end{center}
\caption{Vortex flow around three cylinders with (a) the circulations equal to $0$ on the boundaries of the cylinders, (b) the circulations equal to $-1/3$ on the boundaries of the cylinders.}
\label{fig9}
\end{figure}
\begin{figure}
\begin{center}
\includegraphics[width=0.4\textwidth]{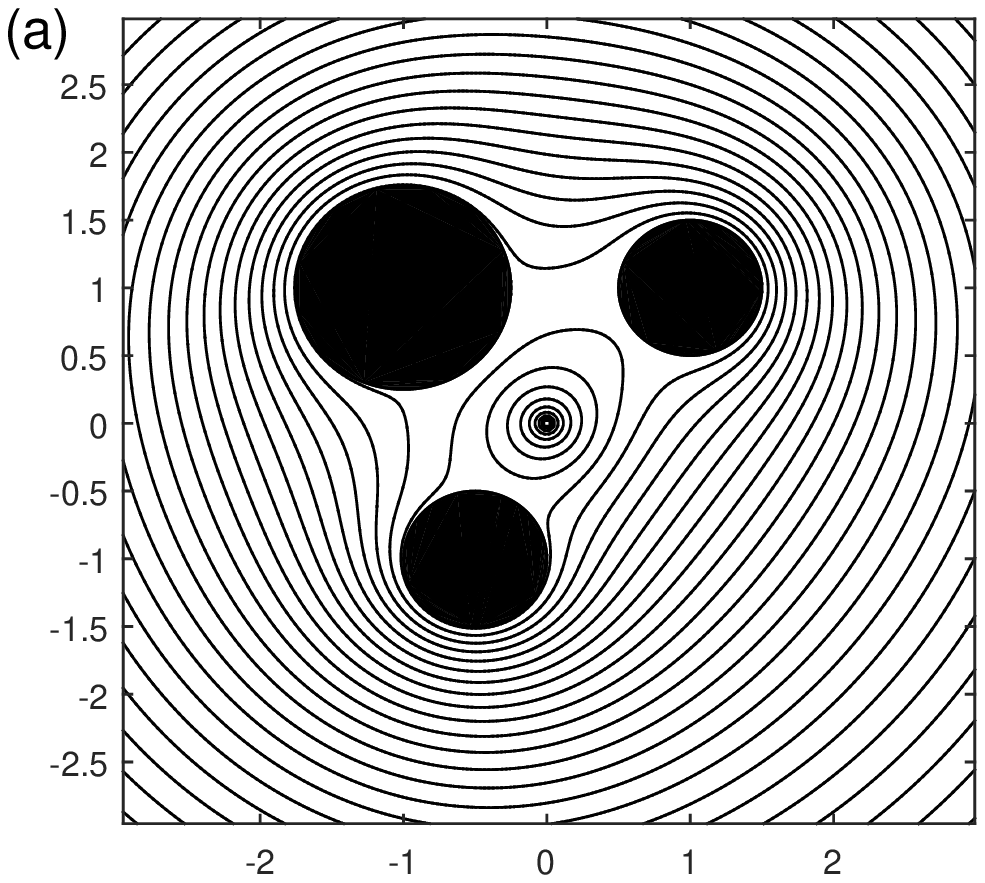}\,\,
\includegraphics[width=0.4\textwidth]{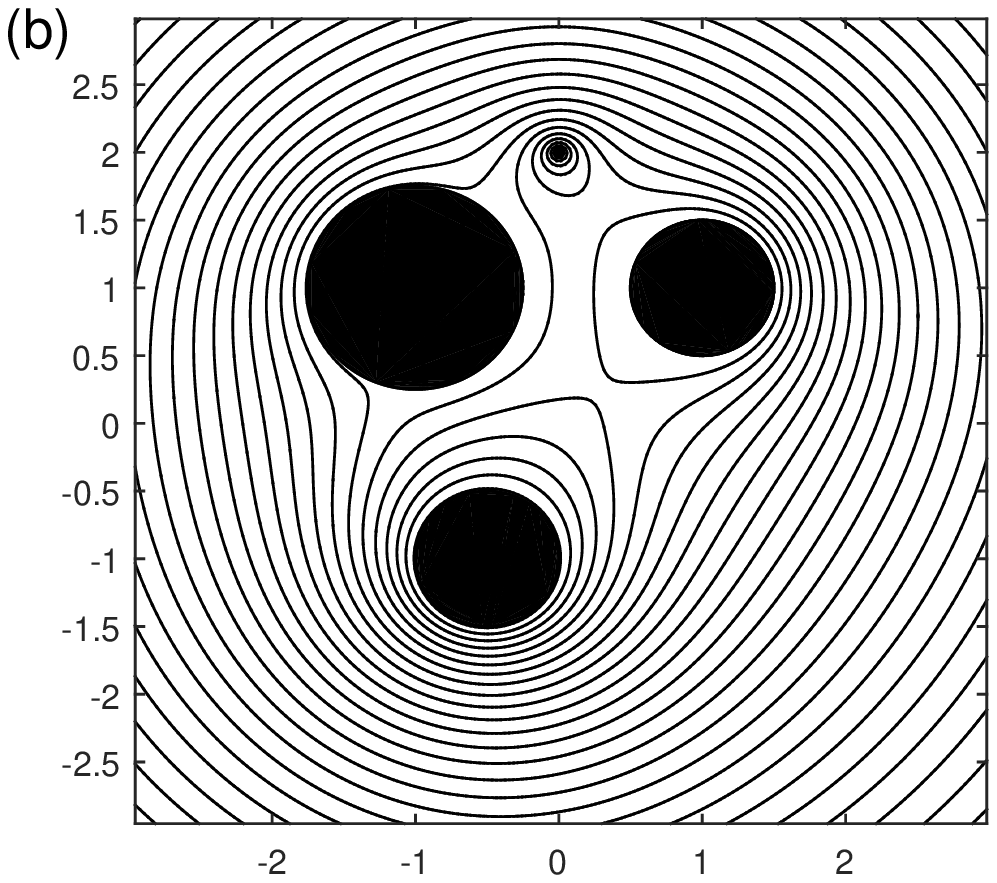}
\end{center}
\caption{Vortex flow around three cylinders with (a) the circulations equal to $0$ on the boundaries of the cylinders and one vortex at the point $z_0=0$, (b) the circulations equal to $0$, $-1$ and $1$ on the boundaries of the cylinders and one vortex at the point $z_0=2i$.}
\label{fig10}
\end{figure}
The graphs of the instantaneous streamlines for the vortex flow around three cylinders are plotted on the figs. \ref{fig9}, \ref{fig10}. The single vortex is located at the point $z_0=0$ and has the circulation $\Gamma_0=1$. On the fig. \ref{fig9} the cylinders have the centers $c_j=e^{ij\pi/3}$, $j=0,1,2$, and the radii $R_j=0.5$. The circulations around all of the cylinders are equal to zero on the fig. \ref{fig9}(a), and are equal to $-1/3$ on the fig. \ref{fig9}(b).  On the fig. \ref{fig10} the cylinders have the centers $c_1=1+i$, $c_2=-1+i$, $c_3=-0.5-i$, and the radii $R_1=0.5$, $R_2=0.75$, and $R_3=0.5$. The circulations around all of the cylinders are equal to zero on the fig. \ref{fig10}(a). On the fig. \ref{fig10}(b) the circulation is equal to $\gamma_1=0$ around the first cylinder, $\gamma_2=-1$ around the second cylinder, and $\gamma_3=1$ around the third cylinder. The vortex with a circulation $\Gamma_0=1$ is located at the point $z_0=0$ for the fig. \ref{fig10}a, and at the point $z_0=2i$ for the fig. \ref{fig10}b.
\begin{figure}
\begin{center}
\includegraphics[width=0.75\textwidth]{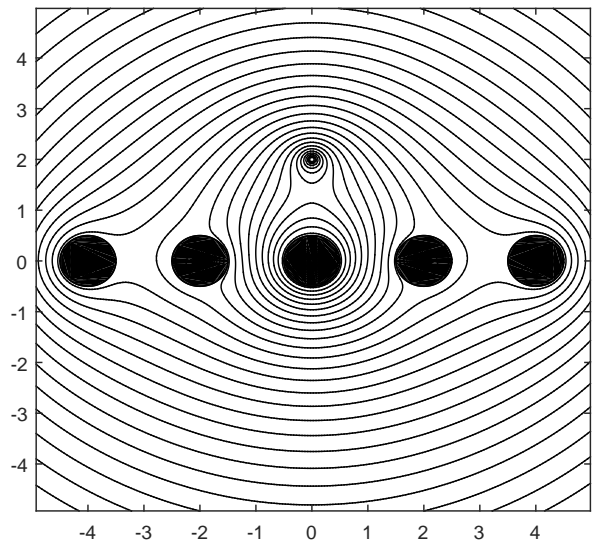}\
\end{center}
\caption{Vortex flow around five cylinders with prescribed circulations.}
\label{fig11}
\end{figure}

The graphs of the instantaneous streamlines for a vortex flow around five cylinders are plotted on the fig. \ref{fig11}. The single vortex is located at the point $z_0=2i$ and has the circulation $\Gamma_0=1$. On the fig. \ref{fig11} the cylinders have the centers $c_1=-4$, $c_2=-2$, $c_3=0$, $c_4=2$, $c_5=4$, and the radii $R_j=0.5$. The circulations are equal to $0$ on the first and the fifth cylinder, equal to $-1$ on the second and the forth cylinder, and equal to $1$ on the third cylinder.

Observe that in all the considered examples, it has been sufficient to consider the symmetry points of the level five at most, thus, the method converges relatively fast. Increasing the number of points beyond this level did not result in a noticeable difference in the pictures of the streamlines of the flow. Comparison of the results for $N=5$ and $N=10$ is given in the table \ref{tab2}. The stream function is computed for the configuration shown on the fig. \ref{fig11}. Observe that taking the symmetry points up to the level $N=5$ already provides us with the first 5 digits after the decimal point.

\begin{center}
\begin{table}[!h]
\begin{tabular}{ |c|c|c| } 
\hline
 Point $z$ & $N=5$ & $N=10$ \\ 
\hline
 $-2-2i$ & $-1.039510891688030$ & $-1.039511060181374$ \\ 
 $4i$ & $-1.127511567288519$ & $-1.127516103881800$ \\ 
 $4-2i$ & $-1.193405902645471$ & $-1.193403567442811$\\
 \hline
\end{tabular}
\label{tab2}
\caption{Comparison of the values of the stream function for $N=5$ and $N=10$.}
\end{table}
\end{center}

\section{Conclusions}

This paper presents a new simple method of study of the vortex generated flows of liquid in domains of arbitrary connectivity using the application of the method of images. It is necessary to observe that the study of the fluid flows in multiply connected domains has a very limited coverage in the scientific literature. To the best of the authors' knowledge, the only other available results are the series of works by Crowdy and Marshall \cite{CrowdyMarshall2005, CrowdyMarshall2005b, CrowdyMarshall2006}. The construction of the stream function presented in this paper is based on taking a limit of a certain functional  sequence which converges to the sought after stream funcion of the fluid flow. The convergence of this functional sequence and its speed are investigated, and the condition of convergence is given as a simple inequality with respect to the geometrical parameters of the flow domain. The limitations of the current study are similar to those in the works of Crowdy and Marshall. The convergence of the presented method is reliable and fast in the case of well-separated cylinders when the connectivity $K$ of the flow domain is not too large. In the cases of very high connectivity $K$ or closely spaced cylinders it may be more efficient to compute the stream function of the flow by using the numerical algorithm proposed by Trefethen \cite{Trefethen2005}.  The results of the current paper can be applied to many practical problems which involve solving the Dirichlet problem for Laplace's equation in multiply connected domains.

\section*{Acknowledgement}
Anna Zemlyanova's research is partially funded through Simon's Foundation Collaboration Grant. This support is gratefully acknowledged.

The authors are grateful to Prof. Hrant Hakobyan for very useful discussions about the nature of the set of the symmetry points.

\end{document}